\documentclass{article}
\usepackage{amsfonts}
\usepackage{hyperref}

\newtheorem{theorem}{Theorem}

\newtheorem{proposition}[theorem]{Proposition}

\newenvironment{proof}[1][Proof]{\noindent\textbf{#1.} }{\ \rule{0.5em}{0.5em}}
\input{tcilatex}

\begin{document}

\title{Hypersurfaces and Codazzi tensors}
\author{Thomas Hasanis and Theodoros Vlachos \\
%EndAName
University of Ioannina\\
E-mail addresses: thasanis@uoi.gr, tvlachos@uoi.gr}
\maketitle

\begin{abstract}
In this paper we deal with the following problem: Let\textit{\ }$%
(M^{n},\left\langle ,\right\rangle )$\textit{\ }be an\textit{\ }$n$\textit{-}%
dimensional Riemannian manifold and\textit{\ }$f:(M^{n},\left\langle
,\right\rangle )\rightarrow \mathbb{R}^{n+1}$\textit{\ }an isometric
immersion. Find all Riemannian metrics on\textit{\ }$M^{n}$\textit{\ }that
can be realized isometrically as immersed hypersurfaces in $\mathbb{R}%
^{n+1}. $ More precisely, given another Riemannian metric\textit{\ }$%
\widetilde{\left\langle ,\right\rangle }$\textit{\ }on\textit{\ }$M^{n},$%
\textit{\ }find necessary and sufficient conditions such that the Riemannian
manifold\textit{\ }$(M^{n},\widetilde{\left\langle ,\right\rangle })$\textit{%
\ }admits an isometric immersion\textit{\ }$\widetilde{f}$\textit{\ }into
the Euclidean space\textit{\ }$\mathbb{R}^{n+1}.$\textit{\ }If such
an\smallskip\ isometric immersion exists, how can one describe\textit{\ }$%
\widetilde{f}$\textit{\ }in terms of\textit{\ }$f?$
\end{abstract}

\section{Introduction and statement of results}

According to a fundamental result due to Nash \cite{Na}, every $n$%
-dimensional Riemannian manifold $(M^{n},\left\langle ,\right\rangle )$\
admits an isometric immersion into the Euclidean space $\mathbb{R}^{n+m},$
for some large $m.$ It is therefore meaningful to ask for isometric
immersions of $(M^{n},\left\langle ,\right\rangle )$ into $\mathbb{R}^{n+m}$
with the lowest possible codimension $m.$ Along this line, Schlaefi in 1873
and later Yau \cite{Y} posed the following conjecture: \textit{any
2-dimensional Riemannian manifold always has a local isometric immersion
into }$\mathbb{R}^{3}.$

It is not always possible for a Riemannian manifold to be realized
isometrically as a hypersurface in a Euclidean space. For example, it is
well known that the hyperbolic space $\mathbb{H}^{n}$ does not admit an
isometric immersion into $\mathbb{R}^{n+1},$ even locally for $n\geq 3.$\ In 
\cite{V},\ Vilms by using the method of bivectors gave necessary and
sufficient conditions for the existence of local isometric immersions of $%
(M^{n},\left\langle ,\right\rangle )$ into the Euclidean space $\mathbb{R}%
^{n+1}.$\ Barbosa, do Carmo and Dajzcer \cite{B1,B2} gave necessary and
sufficient conditions for a Riemannian manifold to be minimally immersed as
a hypersurface in a space form.

The aim of this paper is to begin the study of and call attention to the
following problem: Let\textit{\ }$(M^{n},\left\langle ,\right\rangle )$%
\textit{\ }be an\textit{\ }$n$-dimensional Riemannian manifold and\textit{\ }%
$f:(M^{n},\left\langle ,\right\rangle )\rightarrow \mathbb{R}^{n+1}$\textit{%
\ }an isometric immersion. Find all Riemannian metrics on\textit{\ }$M^{n}$%
\textit{\ }that can be realized isometrically as immersed hypersurfaces in $%
\mathbb{R}^{n+1}.$ More precisely, given another Riemannian metric\textit{\ }%
$\widetilde{\left\langle ,\right\rangle }$\textit{\ }on $M^{n},$\textit{\ }%
find necessary and sufficient conditions such that the Riemannian manifold $%
(M^{n},\widetilde{\left\langle ,\right\rangle })$\textit{\ }admits an
isometric immersion $\widetilde{f}$\textit{\ }into the the Euclidean space $%
\mathbb{R}^{n+1}.$\textit{\ }If such an isometric immersion exists, how can
one describe $\widetilde{f}$\textit{\ }in terms of $f?$

Any Riemannian metric $\widetilde{\left\langle ,\right\rangle }$ on $M^{n}$
determines uniquely a non-singular (1,1)-tensor field\textit{\ }$L$ which is
positive definite, self-adjoint with respect to $\left\langle ,\right\rangle 
$ and satisfies $\widetilde{\left\langle X,Y\right\rangle }=\left\langle
LX,Y\right\rangle ,$ for arbitrary tangent vector fields $X,Y.$ Conversely,
every positive definite (1,1)-tensor field\textit{\ }$L$ which is
self-adjoint with respect to $\left\langle ,\right\rangle $ gives rise to a
new Riemannian metric on $M^{n}.$

In case $L=Id,$ where $Id$\ is the identity map, the above problem reduces
to a rigidity question about $f.$ In fact, according to the classical
Beez-Killing theorem, any other isometric immersion $\widetilde{f}%
:(M^{n},\left\langle ,\right\rangle )\rightarrow \mathbb{R}^{n+1}$ coincides
with $f$ up to an isometry of $\mathbb{R}^{n+1},$ provided that the rank of
the shape operator of $f$ is at least three.

When $L=e^{\varphi }Id,$ where $\varphi $\ is a smooth function on $M^{n},$
the above problem reduces to the study of conformally deformable
hypersurfaces. This case has been studied in details by Cartan \cite{C}. A
modern version of Cartan's results was given recently by Dajczer and Tojeiro 
\cite{DT}.

An important class of Riemannian metrics on $M^{n}$ arises in the case where 
$L=Q^{2},$ $Q$ being a non-singular Codazzi tensor.

In this paper, we give a complete answer to the above problem for this class
of Riemannian metrics.

\begin{theorem}
Let\textit{\ }$f:(M^{n},\left\langle ,\right\rangle )\rightarrow \mathbb{R}%
^{n+1}$\textit{\ be an isometric immersion of a simply connected Riemannian
manifold }$(M^{n},\left\langle ,\right\rangle )$ with shape operator $A$
associated with a unit normal vector field $N.$ Let $\widetilde{\left\langle
,\right\rangle }$ be a new metric\ on $M^{n}$ given by $\widetilde{%
\left\langle X,Y\right\rangle }=\left\langle Q^{2}X,Y\right\rangle ,$ where $%
Q$ is a non-singular Codazzi tensor and $X,Y$ are arbitrary tangent vector
fields. Assume that $\func{rank}A\geq 3.$ Then:

(i) The \textit{Riemannian manifold }$(M^{n},\widetilde{\left\langle
,\right\rangle })$ admits \textit{an isometric immersion into }$\mathbb{R}%
^{n+1}$ if and only if $Q$ commutes with $A.$ If $\widetilde{f}:(M^{n},%
\widetilde{\left\langle ,\right\rangle })\rightarrow \mathbb{R}^{n+1}$ is
such \textit{an isometric immersion, then }$\widetilde{f}$ is rigid with
shape operator $\widetilde{A}=\pm Q^{-1}\circ A.$

(ii) If $Q$ commutes with $A,$ then there exist smooth functions $%
g,h:M^{n}\rightarrow \mathbb{R}$ such that $A(\func{grad}g)=-\func{grad}h$
and $QX=\nabla _{X}\func{grad}g-hAX,$ where $X$ is an arbitrary tangent
vector field and $\nabla $ stands for the Levi-Civit\'{a} connection of $%
(M^{n},\left\langle ,\right\rangle ).$ Moreover, any isometric immersion $%
\widetilde{f}:(M^{n},\widetilde{\left\langle ,\right\rangle })\rightarrow 
\mathbb{R}^{n+1}$ is given by $\widetilde{f}=\tau \circ F,$ where $\tau $ is
an isometry of $\mathbb{R}^{n+1}$ and $F=df(\func{grad}g)+hN.$
\end{theorem}

\section{Preliminaries}

Let $(M^{n},\left\langle ,\right\rangle )$ be an $n$-dimensional Riemannian%
\textit{\ }manifold and $f:M^{n}\rightarrow \mathbb{R}^{n+1}$ an isometric
immersion into the Euclidean space $\mathbb{R}^{n+1}$ with induced bundle $%
f^{\ast }(T\mathbb{R}^{n+1}).$ If $N$ is a unit normal vector field, then
any section $Z$ of the bundle $f^{\ast }(T\mathbb{R}^{n+1})$ decomposes into
a tangent and a normal component, namely%
\[
Z=df(Z_{\top })+hN, 
\]%
where $Z_{\top }$ is a tangent vector field and $h$ is a smooth function on $%
M^{n}.$ Denote by $\overline{\nabla }$ the connection of the induced bundle $%
f^{\ast }(T\mathbb{R}^{n+1})$ arising from the usual connection in $\mathbb{R%
}^{n+1}.\ $For tangents vector fields $X,Y$\ of $M^{n},$\ we have the
Weingarten formula 
\[
\overline{\nabla }_{X}N=-df(AX) 
\]%
and the Gauss formula 
\[
\overline{\nabla }_{X}df(Y)=df(\nabla _{X}Y)+\left\langle AX,Y\right\rangle
N, 
\]%
where $A$ is a self-adjoint (1,1)-tensor field known as the shape operator
associated with $N,$ and $\nabla $ is the Levi-Civit\'{a} connection of $%
\left\langle ,\right\rangle .$

Using the Gauss and Weingarten formulas, one can derive the equations of
Gauss and Codazzi which are, respectively, 
\[
R(X,Y)Z=\left\langle AY,Z\right\rangle AX-\left\langle AX,Z\right\rangle AY, 
\]%
\[
\left( \nabla _{X}A\right) Y=\left( \nabla _{Y}A\right) X, 
\]%
for any tangent vector fields $X,Y,Z,$ where $R$ is the curvature tensor of $%
(M^{n},\left\langle ,\right\rangle )$ and $\nabla _{X}A$ is the covariant
derivative of $A.$

Conversely, the fundamental theorem of hypersurfaces ensures that if there
exists a self-adjoint (1,1)-tensor field $A$ on a simply connected
Riemannian manifold $(M^{n},\left\langle ,\right\rangle )$ that fulfils the
Gauss and the Codazzi equations, then $(M^{n},\left\langle ,\right\rangle )$
admits an isometric immersion into the Euclidean space $\mathbb{R}^{n+1}$
with shape operator $A.$

A self-adjoint (1,1)-tensor field $Q$ on the Riemannian manifold $%
(M^{n},\left\langle ,\right\rangle )$ is said to be a Codazzi tensor if it
satisfies the Codazzi equation $\left( \nabla _{X}Q\right) Y=\left( \nabla
_{Y}Q\right) X$ for arbitrary tangent vector fields $X,Y.$ A family of
Codazzi tensors on $(M^{n},\left\langle ,\right\rangle )$ is given by $%
Q=Id-tA,t\in \mathbb{R}.$ Moreover, if $g,h:M^{n}\rightarrow \mathbb{R}$ are
smooth functions that satisfy $A(\func{grad}g)=-\func{grad}h,$ then the
(1,1)-tensor field $Q$ defined by $QX=\nabla _{X}\func{grad}g-hAX$ is a
Codazzi tensor, where $X$ is an arbitrary tangent vector field. In fact,
this follows by a direct computation and the Gauss equation for $f.$

Consider now a new metric\ $\widetilde{\left\langle ,\right\rangle }$ on $%
M^{n}$ given by $\widetilde{\left\langle X,Y\right\rangle }=\left\langle
Q^{2}X,Y\right\rangle ,$ where $Q$ is a non-singular Codazzi tensor. We
recall the well known relation between the metric $\widetilde{\left\langle
,\right\rangle }$ and the corresponding Levi-Civit\'{a} connection $%
\widetilde{\nabla },$ namely%
\begin{eqnarray*}
2\widetilde{\left\langle \widetilde{\nabla }_{Y}X,Z\right\rangle } &=&X%
\widetilde{\left\langle Y,Z\right\rangle }+Y\widetilde{\left\langle
X,Z\right\rangle }-Z\widetilde{\left\langle X,Y\right\rangle }-\widetilde{%
\left\langle [X,Z],Y\right\rangle } \\
&&-\widetilde{\left\langle [Y,Z],X\right\rangle }-\widetilde{\left\langle
[X,Y],Z\right\rangle }.
\end{eqnarray*}%
In view of the Codazzi equation $\left( \nabla _{X}Q\right) Y=\left( \nabla
_{Y}Q\right) X,$ a direct computation yields%
\begin{equation}
\widetilde{\nabla }_{X}Y=Q^{-1}\left( \nabla _{X}\left( QY\right) \right) . 
\tag{2.1}
\end{equation}%
Moreover, the curvature tensor $\widetilde{R}$ of $(M^{n},\widetilde{%
\left\langle ,\right\rangle })$ is given by

\begin{equation}
\widetilde{R}(X,Y)Z=Q^{-1}\left( R(X,Y)QZ\right)  \tag{2.2}
\end{equation}%
for arbitrary tangent vector fields $X,Y,Z.$

\section{Proof of the result}

For the proof of the main result we need the following auxiliary results.

\begin{proposition}
Let\textit{\ }$f:(M^{n},\left\langle ,\right\rangle )\rightarrow \mathbb{R}%
^{n+1}$\textit{\ be an isometric immersion of a simply connected Riemannian
manifold }$(M^{n},\left\langle ,\right\rangle )$ with shape operator $A$
associated with a unit normal vector field $N.$ Let $\widetilde{\left\langle
,\right\rangle }$ be a new metric\ on $M^{n}$ given by $\widetilde{%
\left\langle X,Y\right\rangle }=\left\langle Q^{2}X,Y\right\rangle ,$ where $%
Q$ is a non-singular Codazzi tensor and $X,Y$ are arbitrary tangent vector
fields. Assume that $\func{rank}A\geq 3.$ Then the \textit{Riemannian
manifold }$(M^{n},\widetilde{\left\langle ,\right\rangle })$ admits \textit{%
an isometric immersion into }$\mathbb{R}^{n+1}$ if and only if $Q$ commutes
with $A.$ If $\widetilde{f}:(M^{n},\widetilde{\left\langle ,\right\rangle }%
)\rightarrow \mathbb{R}^{n+1}$ is such \textit{an isometric immersion, then }%
$\widetilde{f}$ is rigid with shape operator $\widetilde{A}=\pm Q^{-1}\circ
A.$
\end{proposition}

\begin{proof}
We assume that there exists an isometric immersion $\widetilde{f}:(M^{n},%
\widetilde{\left\langle ,\right\rangle })\rightarrow \mathbb{R}^{n+1}$ with
shape operator $\widetilde{A}.$ In view of (2.2), the Gauss equation for $%
\widetilde{f}$ yields 
\begin{eqnarray*}
&&\left\langle Q\circ \widetilde{A}\left( Y\right) ,QZ\right\rangle Q\circ 
\widetilde{A}\left( X\right) -\left\langle Q\circ \widetilde{A}\left(
Y\right) ,QZ\right\rangle Q\circ \widetilde{A}\left( X\right) \\
&=&\left\langle AY,QZ\right\rangle AX-\left\langle AX,QZ\right\rangle AY
\end{eqnarray*}%
for any tangent vector fields $X,Y,Z,$ or equivalently 
\begin{equation}
Q\circ \widetilde{A}\left( X\right) \wedge Q\circ \widetilde{A}\left(
Y\right) =AX\wedge AY,  \tag{3.1}
\end{equation}%
where $\wedge $ stands for the wedge product.

We claim that $\ker A=\ker \widetilde{A}.$ Let $e_{1},...,e_{r}$ be an
orthonormal basis of $\left( \ker A\right) ^{\bot }$\ with respect to $%
\left\langle ,\right\rangle ,$ such that $Ae_{i}=k_{i}e_{i},i=1,...,r,$
where $r=\dim \left( \ker A\right) ^{\bot }\geq 3$ and $\left( \ker A\right)
^{\bot }$ stands for the orthogonal complement of the kernel of $A.$ From
(3.1), we get $AX\wedge e_{i}=0,$ for any $X\in \ker \widetilde{A}$ and $%
i=1,...,r.$ Thus $X\in \ker A$ and $\ker \widetilde{A}\subseteq \ker A.$
Conversely, let $X\in \ker A.$ Then (3.1) yields $Q\circ \widetilde{A}\left(
X\right) \wedge Q\circ \widetilde{A}\left( e_{i}\right) =0,$ for $i=1,...,r.$
Since $Q\circ \widetilde{A}\left( e_{i}\right) \neq 0,$ we obtain $Q\circ 
\widetilde{A}\left( X\right) =\rho _{i}Q\circ \widetilde{A}\left(
e_{i}\right) ,$ for some $\rho _{i},i=1,...,r,$ or equivalently, $\widetilde{%
A}\left( X-\rho _{i}e_{i}\right) =0.$ In view of $\ker \widetilde{A}%
\subseteq \ker A,$ we get $AX=\rho _{i}Ae_{i}$ for $i=1,...,r.$ This implies
that $\rho _{i}=0$ for all $i=1,...,r.$ Thus $Q\circ \widetilde{A}(X)=0$ and 
$X\in \ker \widetilde{A}.$ Consequently $\ker A=\ker \widetilde{A}.$

Let $X\in \left( \func{Ker}A\right) ^{\bot }$ and assume that $Q\circ 
\widetilde{A}\left( X\right) ,AX$ are linearly independent. Since $\dim
\left( \func{Ker}A\right) ^{\bot }\geq 3,$ there exists $Y\in \left( \func{%
Ker}A\right) ^{\bot }$ such that $Q\circ \widetilde{A}\left( X\right) ,AX,AY$
are linearly independent. Then from (3.1) we obtain%
\[
Q\circ \widetilde{A}\left( X\right) \wedge Q\circ \widetilde{A}\left(
X\right) \wedge Q\circ \widetilde{A}\left( Y\right) =Q\circ \widetilde{A}%
\left( X\right) \wedge AX\wedge AY\neq 0, 
\]%
which is contradiction. Thus $Q\circ \widetilde{A}\left( X\right) ,AX$ are
linearly dependent for any $X\in \left( \func{Ker}A\right) ^{\bot }$ and
consequently $Q\circ \widetilde{A}\left( X\right) =a(X)AX,$ where $a(X)$
depends on $X.$

We choose an arbitrary basis $X_{1},...,X_{r}$ of $\left( \func{Ker}A\right)
^{\bot }.$ Then we have $Q\circ \widetilde{A}\left( X_{i}\right)
=a(X_{i})AX_{i},i=1,...,r.$ Since $X_{i}+X_{j}\in \left( \func{Ker}A\right)
^{\bot },$ we have $Q\circ \widetilde{A}\left( X_{i}+X_{j}\right)
=a(X_{i}+X_{j})A\left( X_{i}+X_{j}\right) ,$ for any $i,j=1,...,r,$ and
consequently, 
\[
\Big (a(X_{i}+X_{j})-a(X_{i})\Big )AX_{i}+\Big (a(X_{i}+X_{j})-a(X_{j})\Big )%
AX_{j}=0. 
\]%
This means that $a(X_{i}+X_{j})=a(X_{i})=a(X_{j}),$ for all $i,j=1,...,r.$
Moreover, for any real number $\lambda $, using the linearity of $Q\circ 
\widetilde{A},$ we obtain $Q\circ \widetilde{A}(\lambda X)=\lambda Q\circ 
\widetilde{A}(X)=\lambda a(X)AX.$ Because of $Q\circ \widetilde{A}(\lambda
X)=a(\lambda X)A(\lambda X)=\lambda a(\lambda X)A(X),$ we conclude $%
a(\lambda X)=a(X)$. Therefore, the exists a constant $a$ such that $Q\circ 
\widetilde{A}(X)=aA(X)$ for any $X.$ Appealing to (3.1) we get $a=\pm 1$ and
therefore $\widetilde{A}=\pm Q^{-1}\circ A.$ Since $\widetilde{A}$ is
self-adjoint with respect to $\widetilde{\left\langle ,\right\rangle },$ $%
Q^{2}\circ \widetilde{A}$ is self-adjoint with respect to the metric $%
\left\langle ,\right\rangle ,$ and so we obtain $Q\circ A=A\circ Q.$

Conversely, we assume that $Q$ commutes with $A.$ Then $\widetilde{A}:=\pm
Q^{-1}\circ A$ is self-adjoint with respect to the metric $\widetilde{%
\left\langle ,\right\rangle }.$ Moreover, in view of (2.2), $\widetilde{A}$
satisfies the Gauss equation. In addition, it can be easily seen using (2.1)
that $\widetilde{A}$ is a Codazzi tensor. Then according to the fundamental
theorem of hypersurfaces there exists an isometric immersion $\widetilde{f}%
:(M^{n},\widetilde{\left\langle ,\right\rangle })\rightarrow \mathbb{R}%
^{n+1} $ with shape operator $\widetilde{A}=\pm Q^{-1}\circ A.$ Obviously,
due to the Beez-Killing theorem $\widetilde{f}$ is rigid and the proof is
complete.
\end{proof}

\bigskip

The following result is essentially due to Dajczer and Tojeiro \cite{DT1}.

\begin{proposition}
Let\textit{\ }$f:(M^{n},\left\langle ,\right\rangle )\rightarrow \mathbb{R}%
^{n+1}$\textit{\ be an isometric immersion of a simply connected Riemannian
manifold }$(M^{n},\left\langle ,\right\rangle )$ with shape operator $A$
associated with a unit normal vector field $N.$ Assume that $Q$ is a
non-singular Codazzi tensor which commutes with $A.$ Then there exist:

(i) an immersion $F:M^{n}\rightarrow \mathbb{R}^{n+1}$ such that $dF=df\circ
Q,$ and

(ii) smooth functions $g,h:M^{n}\rightarrow \mathbb{R}$ such that $A(\func{%
grad}g)=-\func{grad}h,$ $QX=\nabla _{X}\func{grad}g-hAX,$ where $X$ is an
arbitrary tangent vector field and $F$ is given by $F=df(\func{grad}g)+hN.$

Furthermore, $F$ is an isometric immersion of the \textit{Riemannian
manifold (}$M^{n},\widetilde{\left\langle ,\right\rangle })$ into $\mathbb{R}%
^{n+1},$\ where the metric $\widetilde{\left\langle ,\right\rangle }$\ on $%
M^{n}$ is given by $\widetilde{\left\langle X,Y\right\rangle }=\left\langle
Q^{2}X,Y\right\rangle ,$ $X,Y$ being arbitrary tangent vector fields.
\end{proposition}

\begin{proof}
(i) Let $e_{1},...,e_{n+1}$ be the standard orthonormal basis of $\mathbb{R}%
^{n+1}.\mathbb{\ }$We consider the 1-forms $\omega _{i}:=\left\langle
df\circ Q,e_{i}\right\rangle ,i=1,...,n+1.$\ Then for any\ tangent vector
fields $X,Y,$ we get 
\begin{eqnarray*}
d\omega _{i}(X,Y) &=&X\left( \omega _{i}(Y)\right) -Y\left( \omega
_{i}(X)\right) -\omega _{i}\left( [X,Y]\right) \\
&=&\left\langle \overline{\nabla }_{X}df(QY),e_{i}\right\rangle
-\left\langle \overline{\nabla }_{Y}df(QX),e_{i}\right\rangle -\left\langle
df\left( Q[X,Y]\right) ,e_{i}\right\rangle .
\end{eqnarray*}%
Using the Gauss formula, we obtain%
\begin{eqnarray*}
d\omega _{i}(X,Y) &=&\left\langle df\left( \nabla _{X}(QY)-\nabla
_{Y}(QX)-Q([X,Y])\right) ,e_{i}\right\rangle \\
&&+\left\langle [Q,A]X,Y\right\rangle \left\langle N,e_{i}\right\rangle ,
\end{eqnarray*}%
for each$\ i=1,...,n+1,$ where $[Q,A]=Q\circ A-A\circ Q.$ Since $Q$ is a
Codazzi tensor and commutes with $A,$ we deduce that $\omega _{i}$'s are
closed. Thus there exist functions $u_{i}:M^{n}\rightarrow \mathbb{R}$ such
that $\omega _{i}=du_{i}\ i=1,...,n+1.$ Then the map $F:M^{n}\rightarrow 
\mathbb{R}^{n+1}$ given by $F:=\sum_{i=1}^{n+1}u_{i}e_{i}$\ clearly
satisfies $dF=df\circ Q.$

(ii) We view the immersion $F$ introduced in (i) as a section of the induced
bundle $f^{\ast }(T\mathbb{R}^{n+1}).$ Decomposing into a tangent and a
normal component, we get%
\begin{equation}
F=df(Z)+hN,  \tag{3.2}
\end{equation}%
where $Z$ is a tangent vector field and $h$ is a smooth function on $M^{n}.$
Then on account of Gauss and Weingarten formulas, for any tangent vector
field $X,$ we have%
\[
dF(X)=df\left( \nabla _{X}Z-hAX\right) +\left( \left\langle
AX,Z\right\rangle +Xh\right) N, 
\]%
or equivalently 
\[
df\left( QX-\nabla _{X}Z+hAX\right) =\left( \left\langle AX,Z\right\rangle
+Xh\right) N. 
\]%
From this we find that $QX=\nabla _{X}Z-hAX$ and $AZ=-\func{grad}h.$
Furthermore, we deduce that the 1-form $\omega $ defined by $\omega
(X)=\left\langle X,Z\right\rangle $ is closed. Since $M^{n}$ is simply
connected, there exists a function $g:M^{n}\rightarrow \mathbb{R}$ such that 
$\omega =dg,$ or equivalently $Z=\func{grad}g.$ Then (3.2) is written as $%
F=df(\func{grad}g)+hN.$

Moreover, the metric $\widetilde{\left\langle ,\right\rangle }$ induced on $%
M^{n}$ is given by 
\[
\widetilde{\left\langle X,Y\right\rangle }=\left\langle
dF(X),dF(Y)\right\rangle =\left\langle df\left( QX\right) ,df\left(
QY\right) \right\rangle =\left\langle Q^{2}X,Y\right\rangle , 
\]%
for any\ tangent vector fields $X,Y.$ This completes the proof.
\end{proof}

\bigskip

\begin{proof}[Proof of the Theorem 1]
Let\textit{\ }$f:(M^{n},\left\langle ,\right\rangle )\rightarrow \mathbb{R}%
^{n+1}$\textit{\ }be an isometric immersion of the simply connected
Riemannian manifold $(M^{n},\left\langle ,\right\rangle )$ with shape
operator $A$ associated with a unit normal vector field $N.$ \smallskip We
consider the Riemannian metric\ $\widetilde{\left\langle ,\right\rangle }$
on $M^{n}$ given by $\widetilde{\left\langle X,Y\right\rangle }=\left\langle
Q^{2}X,Y\right\rangle ,$ where $Q$ is a non-singular Codazzi tensor and $X,Y$
are arbitrary tangent vector fields.

(i) According to Proposition 2, the Riemannian manifold $(M^{n},\widetilde{%
\left\langle ,\right\rangle })$ admits an isometric immersion into\textit{\ }%
$\mathbb{R}^{n+1}$ if and only if $Q$ commutes with $A.$ Moreover, such an
isometric immersion $\widetilde{f}:(M^{n},\widetilde{\left\langle
,\right\rangle })\rightarrow \mathbb{R}^{n+1}$ is rigid with shape operator $%
\widetilde{A}=\pm Q^{-1}\circ A.$

(ii) If $Q$ commutes with $A,$ then Proposition 3 implies that there exists
an isometric immersion $F:(M^{n},\widetilde{\left\langle ,\right\rangle }%
)\rightarrow \mathbb{R}^{n+1}.$ Moreover, there exist smooth functions $%
g,h:M^{n}\rightarrow \mathbb{R}$ such that $A(\func{grad}g)=-\func{grad}h,$ $%
QX=\nabla _{X}\func{grad}g-hAX,$ where $X$ is an arbitrary tangent vector
field, and $F=df(\func{grad}g)+hN.$ In addition, since $\func{rank}A\geq 3,$
by the Beez-Killing theorem, any isometric immersion $\widetilde{f}:(M^{n},%
\widetilde{\left\langle ,\right\rangle })\rightarrow \mathbb{R}^{n+1}$ is
given by $\widetilde{f}=\tau \circ F,$ where $\tau $ is an isometry of $%
\mathbb{R}^{n+1}$.
\end{proof}

\bigskip \textbf{Remark 4 }It is worth noticing that $f$ and $F$ in
Proposition 3 have\smallskip\ the same Gauss map, and consequently the
immersions $f$ and $\widetilde{f}$ in Theorem 1 have congruent Gauss maps.

\bigskip

The following two examples illustrate the main result.

\bigskip

\textbf{Example 5 }Let\textit{\ }$f:(M^{n},\left\langle ,\right\rangle
)\rightarrow \mathbb{R}^{n+1}$\textit{\ }be an isometric immersion of a
simply connected Riemannian manifold $(M^{n},\left\langle ,\right\rangle )$
with shape operator $A$ corresponding to a unit normal vector field $N.$ We
consider the tensor field $Q:=Id-tA,t\in \mathbb{R},$ where $t$ is chosen
such that $Q$ is non-singular. Obviously, $Q$ is a Codazzi tensor that
commutes with $A.$ According to our results, the Riemannian manifold $(M^{n},%
\widetilde{\left\langle ,\right\rangle })$ is realized isometrically as an
immersed hypersurface in $\mathbb{R}^{n+1},$ where the metric $\widetilde{%
\left\langle ,\right\rangle }$ is given by $\widetilde{\left\langle
X,Y\right\rangle }=\left\langle Q^{2}X,Y\right\rangle ,$ and $X,Y$ are
arbitrary tangent vector fields.

Viewing $f$ as a section of the induced bundle $f^{\ast }(T\mathbb{R}%
^{n+1}), $ and decomposing into a tangent and a normal component, we get%
\begin{equation}
f=df(x_{\top })+pN,  \tag{3.3}
\end{equation}%
where $x_{\top }$ is a tangent vector field and $p=\left\langle
f,N\right\rangle $ is the support function of $f.$ Differentiating (3.3)
with respect to a tangent vector field $X$ and using Gauss and Weingarten
formulas, we obtain 
\[
\nabla _{X}x_{\top }=X+pAX\ \ \text{and}\ \ Ax_{\top }=-\func{grad}p. 
\]%
Moreover we, easily, see that $\func{grad}(\frac{1}{2}\left\vert
f\right\vert ^{2})=x_{\top }.$ Consequently the functions $\frac{1}{2}%
\left\vert f\right\vert ^{2}$ and $p$ satisfy 
\[
A\func{grad}(\frac{1}{2}\left\vert f\right\vert ^{2})=-\func{grad}p\ \ \text{%
and}\ \ QX=\nabla _{X}\func{grad}(\frac{1}{2}\left\vert f\right\vert
^{2})-\left( p+t\right) AX. 
\]%
So the functions $g:=\frac{1}{2}\left\vert f\right\vert ^{2}$ and $h:=p+t$
can be used to construct the isometric immersion $F:(M^{n},\widetilde{%
\left\langle ,\right\rangle })\rightarrow \mathbb{R}^{n+1},$ that is $F=df(%
\func{grad}g)+hN.$ In view of (3.3), it follows that $F=f+tN,$ a parallel
hypersurface to $f.$

\bigskip

\textbf{Example 6 }Let\textit{\ }$f:(M^{n},\left\langle ,\right\rangle
)\rightarrow \mathbb{R}^{n+1}$\textit{\ }be an isometric immersion of a
simply connected Riemannian manifold $(M^{n},\left\langle ,\right\rangle )$
with shape operator $A$ associated with a unit normal vector field $N.$
Consider a constant vector $a$ in $\mathbb{R}^{n+1}\ $and view it as a
section of the induced bundle $f^{\ast }(T\mathbb{R}^{n+1}).$ We decompose $%
a $ as 
\begin{equation}
a=df(a_{\top })+\left\langle N,a\right\rangle N.  \tag{3.4}
\end{equation}%
Differentiating with respect to a tangent vector field $X$ and using Gauss
and Weingarten formulas, we obtain 
\[
\nabla _{X}a_{\top }=\left\langle N,a\right\rangle AX\ \ \text{and}\ \
Aa_{\top }=-\func{grad}\left\langle N,a\right\rangle . 
\]%
Since $\func{grad}\left\langle f,a\right\rangle =a_{\top },$ we have$\ A%
\func{grad}\left\langle f,a\right\rangle =-\func{grad}\left\langle
N,a\right\rangle .$ So the functions $g:=\left\langle f,a\right\rangle $ and 
$h:=\left\langle N,a\right\rangle +1$ satisfy 
\[
A\func{grad}g=-\func{grad}h\ \ \text{and}\ \ QX=\nabla _{X}\func{grad}%
g-\left( \left\langle N,a\right\rangle +1\right) AX=-AX, 
\]%
and can be used to construct the isometric immersion $F:(M^{n},\widetilde{%
\left\langle ,\right\rangle })\rightarrow \mathbb{R}^{n+1}$ where $%
\widetilde{\left\langle X,Y\right\rangle }=\left\langle
A^{2}X,Y\right\rangle ,$ and $X,Y$ are arbitrary tangent vector fields. In
this case, we have $F=df(\func{grad}g)+hN=df(a_{\top })+\left( \left\langle
N,a\right\rangle +1\right) N=a+N,$ the Gauss map of $f$ followed by a
parallel translation.

\bigskip

In the following remark we discuss the uniqueness of the pair of functions $%
g,h,$ for a given Codazzi tensor $Q$ which commutes with the shape operator $%
A.$

\bigskip

\textbf{Remark 7 }Under the assumptions of Theorem 1, we suppose that there
exist two pairs of functions $\left( g,h\right) ,\left( g_{1},h_{1}\right) ,$
such that 
\[
A(\func{grad}g)=-\func{grad}h,A(\func{grad}g_{1})=-\func{grad}h_{1}
\]%
and 
\[
QX=\nabla _{X}\func{grad}g-hAX=\nabla _{X}\func{grad}g_{1}-h_{1}AX,
\]%
where $X$ is an arbitrary tangent vector field. Then, according to
Proposition 3, the immersions $F=df(\func{grad}g)+hN$ and $F_{1}=df(\func{%
grad}g_{1})+h_{1}N$ induce the same metric on $M^{n}$ and satisfy $%
dF=df\circ Q=dF_{1}.$ Thus $F_{1}=F+a,$ for a constant vector $a.$ Viewing $a
$ as a section of the induced bundle $f^{\ast }(T\mathbb{R}^{n+1}),$ we
decompose as in (3.4). Since $F=df(\func{grad}g)+hN$ and $F_{1}=df(\func{grad%
}g_{1})+h_{1}N,$ from (3.3) and $F_{1}=F+a,$ we obtain $g_{1}=g+\left\langle
f,a\right\rangle +c$ and $h_{1}=h+\left\langle N,a\right\rangle ,$ where $c$
is a real constant. Thus the pair of functions $g,h$ are uniquely determined
up to the functions $\left\langle f,a\right\rangle $ and $\left\langle
N,a\right\rangle .$

%TCIMACRO{%
%\TeXButton{Undefine the chapter macro (see explanation below)}{\let\chapter\undefined}}%
%BeginExpansion
\let\chapter\undefined%
%EndExpansion
\appendix{}

\end{document}